# On the Existence of Periodic Solutions with Applications to Extremum-Seeking

**Iasson Karafyllis[*] and Miroslav Krstic[**]**

[*]Dept. of Mathematics, National Technical University of Athens, Zografou Campus, 15780, Athens, Greece,
email: iasonkar@central.ntua.gr; iasonkaraf@gmail.com

[**]Dept. of Mechanical and Aerospace Eng., University of California, San Diego, La Jolla, CA 92093-0411, U.S.A., email: krstic@ucsd.edu

**Abstract**

This paper provides two results that can be useful in the study of the existence and the stability properties of a periodic solution for a given dynamical system. The first result deals with scalar time-periodic systems and establishes the equivalence of the existence of a periodic solution and the existence of a bounded solution. The second result provides sufficient conditions for the existence and the stability of a periodic solution for a time-periodic dynamical system. Both results are applied to extremum seeking problems for a static output map with no plant dynamics and novel non-local results are provided without the use of averaging theorems and singular perturbation arguments.

## 1. Introduction

Periodic solutions are ubiquitous in nature as well as in dynamical systems. However, the study of the existence and the stability properties of a periodic solution for a given dynamical system is far more difficult than the corresponding study for equilibrium points or sets of equilibria (see for instance [3, 4]).

In control systems theory, periodic solutions and time-periodic dynamical systems arise naturally in extremum seeking control. Extremum seeking control is a topic that has attracted the attention of many researchers in control theory from the 50s (see for instance [10, 11, 12, 14]) and has recently been extended to various directions (see [2, 8, 13, 16]). The study of extremum seeking control is based mainly on averaging theorems (see [4, 5]) and singular perturbation results (see [5]) which allow the derivation of local (or semi-global as in [18]) results for the existence of an attracting periodic solution.

The present paper provides two main results that can be useful in the study of the existence and the stability properties of a periodic solution for a given dynamical system:

1) Theorem 1 deals with scalar time-periodic systems and establishes the equivalence of the existence of a periodic solution and the existence of a bounded solution,
2) Theorem 2 deals with time-periodic systems and provides sufficient conditions for the existence and the stability of a periodic solution.

Both main results have some important consequences. First of all, Theorem 1 can be used in order to show that a periodic orbit is contained in the omega limit set for a time-invariant dynamical system (see Section 3 below). Secondly, Theorem 2 is related to the recent results in [9, 15] that provide sufficient conditions for the entrainment to periodic external inputs (see the discussion below Theorem 2 below). The existence of periodic solutions of time-invariant control systems that are subject to periodic external inputs plays a crucial role in biology (as argued in [9, 15]) and is also related to stability notions for control systems (see for instance the proof in [17] of the existence of a periodic solution for an Input-to-State Stable control system under a periodic input).

However, both Theorem 1 and Theorem 2 can also be applied in a direct way to extremum seeking problems for a static output map with no plant dynamics: this is the topic of Section 4 below, where novel results are provided (Theorem 4 and Proposition 1) without the use of averaging theorems and singular perturbation results. The obtained results are not local and clarify completely the nature of the dynamics of extremum seeking systems. We are not aware of similar results in the literature.

The structure of the paper is as follows. Section 2 provides our main results (Theorem 1 and Theorem 2) for time-periodic systems. Section 3 is devoted to the application of the main results to time-invariant systems while Section 4 deals with the application of the main results to extremum seeking problems for a static output map with no plant dynamics. All proofs are provided in Section 5. Finally, some concluding remarks are given in Section 6.

**Notation.** Throughout this paper, we adopt the following notation.

* $\mathbb{R}_+ := [0, +\infty)$.

* Let $D \subseteq \mathbb{R}^n$ be a non-empty open set and let $S \subseteq \mathbb{R}^n$ be a set that satisfies $D \subseteq S \subseteq \bar{D}$. By $C^0(S; \Omega)$, we denote the class of continuous functions on $S$, which take values in $\Omega \subseteq \mathbb{R}^m$. By $C^k(S; \Omega)$, where $k \geq 1$ is an integer, we denote the class of functions on $S \subseteq \mathbb{R}^n$, which take values in $\Omega \subseteq \mathbb{R}^m$ and have continuous derivatives of order $k$. In other words, the functions of class $C^k(S; \Omega)$ are the functions which have continuous derivatives of order $k$ in $D$ that can be continued continuously to all points in $S$. When $\Omega = \mathbb{R}$ then we write $C^0(S)$ or $C^k(S)$.

* For a time invariant dynamical system of the form $\dot{x} = f(x)$, $x \in \mathbb{R}^n$ where $f : \mathbb{R}^n \to \mathbb{R}^n$ is a locally Lipschitz vector field we can define the omega limit set of $x_0 \in \mathbb{R}^n$. Let $x(t)$ be the unique solution at time $t \geq 0$ of the initial-value problem $\dot{x} = f(x)$ with initial condition $x(0) = x_0 \in \mathbb{R}^n$ and suppose that is defined for all $t \geq 0$. We define $\omega(x_0) \subseteq \mathbb{R}^n$ to be the $\omega-$limit set of $x_0 \in \mathbb{R}^n$ using the formula:

$$\omega(x_0) = \left\{ y \in \mathbb{R}^n : \exists \{t_k \geq 0 : k = 1, 2, ...\} \text{ with } \lim_{k \to +\infty}(t_k) = +\infty \text{ and } \lim_{k \to +\infty}(x(t_k)) = y \right\} \quad (1.1)$$

* Let $|\cdot|$ denote a norm of $\mathbb{R}^n$. The induced matrix norm $\|\cdot\|$ is defined for a matrix $G \in \mathbb{R}^{n \times n}$ by the formula $\|G\| := \sup\{|Gx| : |x| \leq 1\}$. The associated logarithmic norm $\mu(G)$ of $G \in \mathbb{R}^{n \times n}$ is defined by $\mu(G) := \lim_{h \to 0^+}\left(\dfrac{\|I_n + hG\| - 1}{h}\right)$, where $I_n$ is the unit $n \times n$ matrix.

* For $x \in \mathbb{R}$ we define the integer part of $x$ denoted by $\lfloor x \rfloor$ as the largest integer that is less than or equal to $x$, i.e., $\lfloor x \rfloor = \max\{q \in \mathbb{Z} : q \leq x\}$.

* By $KL$ we denote the set of all continuous functions $\sigma : \mathbb{R}_+ \times \mathbb{R}_+ \to \mathbb{R}_+$ with the following properties: (i) for each $t \geq 0$ the function $a(s) = \sigma(s,t)$ is an increasing function with $a(0) = 0$; (ii) for each $s \geq 0$, the function $\sigma(s,\cdot)$ is non-increasing with $\lim_{t \to +\infty}(\sigma(s,t)) = 0$.

## 2. Main Results

Let $f \in C^1(\mathbb{R}_+ \times \mathbb{R}^n)$ be given. Consider the time-varying system given by the differential equation

$$\dot{x}(t) = f(t, x(t)) \tag{2.1}$$

where $x(t) \in \mathbb{R}^n$ is the state variable.

Let $t_0 \geq 0$, $x_0 \in \mathbb{R}^n$ be given. By $\varphi(t, t_0, x_0)$ we denote the unique solution of (2.1) at time $t \geq t_0$ with initial condition $x(t_0) = x_0$. For each $t_0 \geq 0$, $x_0 \in \mathbb{R}^n$ there exists $t_{\max} \in (t_0, +\infty]$ for which the solution $\varphi(t, t_0, x_0)$ is defined for $t \in [t_0, t_{\max})$. Moreover, if $t_{\max} < +\infty$ then $\limsup_{t \to t_{\max}^-}(|\varphi(t, t_0, x_0)|) = +\infty$, where $|\cdot|$ denotes a norm of $\mathbb{R}^n$. The map $\varphi$ satisfies $\varphi(t_0, t_0, x_0) = x_0$ (identity property) and $\varphi(t, t_0, x_0) = \varphi(t, \tau, \varphi(\tau, t_0, x_0))$ (semigroup property) for all $t_0 \geq 0$, $x_0 \in \mathbb{R}^n$ and $t_0 \leq \tau \leq t < t_{\max}$.

We next assume that there exists $T > 0$ such that the following identity holds for all $t \geq 0$, $x \in \mathbb{R}^n$:

$$f(t + T, x) = f(t, x) \tag{2.2}$$

Identity (2.2) means that system (2.1) is $T$–periodic and the following fact is valid as long as all involved quantities are well defined:

**Fact I:** $\varphi(t + kT, t_0 + kT, x_0) = \varphi(t, t_0, x_0)$, for all $k \in \mathbb{N}$.

When $n = 1$, we obtain for all $t_0 \geq 0$, $x_0 \in \mathbb{R}$ and $t \in [t_0, t_{\max})$:

$$\frac{\partial \varphi}{\partial t}(t,t_0,x_0) = f\left(t,\varphi(t,t_0,x_0)\right)$$

$$\frac{\partial \varphi}{\partial x_0}(t,t_0,x_0) = \exp\left(\int_{t_0}^{t} \frac{\partial f}{\partial x}\left(s,\varphi(s,t_0,x_0)\right)ds\right) \quad (2.3)$$

Since (2.3) gives that $\frac{\partial \varphi}{\partial x_0}(t,t_0,x_0) > 0$, we conclude that the following fact is valid as long as all involved quantities are well defined:

**Fact II:** For $n=1$, $x, y \in \mathbb{R}$: $x < y \Rightarrow \varphi(t,t_0,x) < \varphi(t,t_0,y)$

Using Fact I and Fact II we obtain the following result for system (2.1).

**Theorem 1:** *Let $f \in C^1(\mathbb{R}_+ \times \mathbb{R})$ be a map that satisfies (2.2) for all $t \geq 0$, $x \in \mathbb{R}$. Consider system (2.1) with $n=1$. Let $t_0 \geq 0$, $x_0 \in \mathbb{R}$ be given and suppose that the solution $\varphi(t,t_0,x_0)$ of (2.1) exists for all $t \geq t_0$ and is bounded on $[t_0,+\infty)$. Then there exists a $T$-periodic solution $x^*(t)$ of (2.1) such that $\lim_{t \to +\infty}\left(\left|\varphi(t,t_0,x_0) - x^*(t)\right|\right) = 0$.*

Theorem 1 can be used in three ways:

1) the existence of a periodic solution of (2.1) can be proved by showing that there exists a bounded solution for (2.1),

2) the fact that all bounded solutions of (2.1) converge as $t \to +\infty$ to some periodic solution (but not necessarily the same periodic solution; you may have more than one periodic solutions) can be established. Thus, the most complicated limit set for a scalar ODE with time-periodic right-hand side is a periodic solution, and

3) it can be shown that if a solution of (2.1) does not converge to a periodic solution then the solution is necessarily unbounded (possibly with a finite escape time).

Theorem 1 is proved in Section 5 with the use of the so-called return (Poincaré) maps (see [3]).

The following example shows how easily Theorem 1 can be applied to certain planar systems and guarantee the existence of non-trivial periodic solutions (not necessarily limit cycles).

**Example 1:** Consider the planar system
$$\dot{x}_1 = f(x)x_1 - \omega(x)x_2$$
$$\dot{x}_2 = \beta^2 \omega(x)x_1 + f(x)x_2 \quad (2.4)$$

where $\omega, f : \mathbb{R}^2 \to \mathbb{R}$ are $C^1$ functions with $\omega(x) \neq 0$ for all $x \in \mathbb{R}^2$ and $\beta > 0$ is a constant. System (2.4) in the elliptic coordinates

$$x_1 = r\cos(\theta)$$
$$x_2 = \beta r \sin(\theta) \tag{2.5}$$

is described by the differential equations

$$\dot{r} = r f(r\cos(\theta), \beta r \sin(\theta))$$
$$\dot{\theta} = \beta \omega(r\cos(\theta), \beta r \sin(\theta)) \tag{2.6}$$

Since $\dot{\theta} \neq 0$ (a consequence of (2.6) and the fact that $\omega(x) \neq 0$ for all $x \in \mathbb{R}^2$) we can look for solutions of (2.6) of the form $r(t) = r(\theta(t))$ which satisfy the differential equation:

$$\frac{dr}{d\theta} = \frac{r f(r\cos(\theta), \beta r \sin(\theta))}{\beta \omega(r\cos(\theta), \beta r \sin(\theta))} \tag{2.7}$$

If we also want to exclude the trivial periodic solution $x=0$ (or $r=0$) then we can apply the transformation $r = e^z$ with $z \in \mathbb{R}$ that allows us to obtain from (2.7) the differential equation

$$\frac{dz}{d\theta} = F(\theta, z) := \frac{f(e^z \cos(\theta), \beta e^z \sin(\theta))}{\beta \omega(e^z \cos(\theta), \beta e^z \sin(\theta))} \tag{2.8}$$

Clearly, the function $F(\theta, z)$ defined by (2.8) is a $C^1$ function which is $2\pi$-periodic in $\theta$. A periodic solution $z = \varphi(\theta)$ of (2.8) directly gives a periodic solution of (2.4). This is done by substituting in the equation for $\theta$ and solving the differential equation $\dot{\theta}(t) = \beta \omega\big(\exp(\varphi(\theta(t)))\cos(\theta(t)), \beta \exp(\varphi(\theta(t)))\sin(\theta(t))\big)$ with initial condition $\theta(0) = 0$. When $\omega(x) > 0$ for all $x \in \mathbb{R}^2$ then there exists $T > 0$ such that $\theta(T) = 2\pi$ and when $\omega(x) < 0$ for all $x \in \mathbb{R}^2$ then there exists $T > 0$ such that $\theta(T) = -2\pi$. It is straightforward to verify that the solution $x_1(t) = \exp(\varphi(\theta(t)))\cos(\theta(t))$, $x_2(t) = \beta \exp(\varphi(\theta(t)))\sin(\theta(t))$ for $t \in [0, T]$ is a periodic solution of period $T > 0$ of system (2.4).

In order to obtain a $2\pi$-periodic solution of (2.8) we can use Theorem 1. By showing that there exists a bounded solution of (2.8), Theorem 1 implies the existence of a $2\pi$-periodic solution of (2.8). Consequently, in order to show the existence of a periodic solution of (2.4) it suffices to show the existence of a bounded solution of (2.8). This can be done very easily in some cases: for example, when $\omega(x) > 0$ for all $x \in \mathbb{R}^2$ and there exists $R > 0$ such that $f(r\cos(\theta), \beta r \sin(\theta)) \leq 0$ for all $r \geq R$. Moreover, if (2.8) has no equilibrium point then we can be sure that the obtained periodic solution is non-trivial (it is not an equilibrium point). ◁

Our second main result is a result that guarantees the existence of a periodic solution that attracts all solutions from a given set of initial conditions.

**Theorem 2:** *Let $f \in C^1\left(\mathbb{R}_+ \times \mathbb{R}^n\right)$ be a map that satisfies (2.2) for all $t \geq 0$, $x \in \mathbb{R}^n$. Suppose that there exists a set $S \subseteq \mathbb{R}^n$ and a convex, compact set $C \subseteq \mathbb{R}^n$ such that for all $t_0 \geq 0$, $x_0 \in S$, there exists $\tau \geq t_0$ with $\varphi(t, t_0, x_0) \in C$ for all $t \geq \tau$. Moreover, suppose that there exists a $T$-periodic function $p \in C^0(\mathbb{R}_+)$ with $\int_0^T p(s)ds < 0$ such that $\mu\left(\frac{\partial f}{\partial z}(t, z)\right) \leq p(t)$ for all $t \geq 0$, $z \in C$, where $\mu(\cdot)$ denotes the logarithmic norm associated with a norm $|\cdot|$ of $\mathbb{R}^n$. Then there exists a unique $T$-periodic solution $x^*(t)$ of (2.1) with $x^*(t) \in C$ for all $t \in [0, T]$ such that $\lim_{t \to +\infty} \left(\left|\varphi(t, t_0, x_0) - x^*(t)\right|\right) = 0$ for all $t_0 \geq 0$, $x_0 \in S$.*

It should be noted for system (2.1) with $n = 1$, Fact II has important consequences for Theorem 2: (i) the fact that the set $C \subseteq \mathbb{R}$ has to be convex and compact implies that the set $C$ is an interval of the form $[a, b]$, (ii) due to Fact II and the fact that $C = [a, b]$ it also follows that the set $S \subseteq \mathbb{R}$ in Theorem 2 must also be a connected set, i.e., $S \subseteq \mathbb{R}$ is necessarily an interval. It should also be noted that for system (2.1) with $n = 1$ we have $\mu\left(\frac{\partial f}{\partial z}(t, z)\right) = \frac{\partial f}{\partial z}(t, z)$.

The proof of Theorem 2 is provided in Section 5 and is again based on the estimation of the return (Poincaré) map for system (2.1).

A few remarks are needed here for Theorem 2. Theorem 2 can be compared with Theorem 2 and Theorem 9 in [15]. Here, instead of assuming that (2.1) is forward complete on a closed, positively invariant, convex (or K-reachable) set $C \subseteq \mathbb{R}^n$ and $f$ is infinitesimally contracting with contraction rate $c^2 > 0$ on $C$, i.e., $\mu\left(\frac{\partial f}{\partial z}(t, z)\right) \leq -c^2 < 0$ for all $t \geq 0$, $z \in C$, we assume that $\mu\left(\frac{\partial f}{\partial z}(t, z)\right) \leq p(t)$ for all $t \geq 0$, $z \in C$. Notice that the function $p(t)$ is allowed to take positive values since the only requirement for $p(t)$ is $\int_0^T p(s)ds < 0$. Furthermore, here we do not assume that the set $C \subseteq \mathbb{R}^n$ is a positively invariant set: we only assume that $C \subseteq \mathbb{R}^n$ "traps" the solutions initiated from the set $S \subseteq \mathbb{R}^n$. On the other hand, here we assume that the set $C \subseteq \mathbb{R}^n$ is compact, and we do not assume forward completeness. Theorem 2 can also be compared with Proposition 1 in [9] which assumes the existence of a compact, positively invariant and convex set $C \subseteq \mathbb{R}^n$ and that system (2.1) on $C$ is contractive after a small overshoot and short transient (see [9] for a precise definition of a contractive after a small overshoot and short transient system). It should be noted here that the assumptions of Theorem 2 do not guarantee that system (2.1) on $C$ is contractive after a small overshoot and short transient.

It should also be noticed that due to all the facts that we mentioned above, the extremum seeking application of Theorem 2 that is presented in Section 4 below cannot be handled by Theorem 2 and Theorem 9 in [15] or Proposition 1 in [9]. Theorem 2 takes advantage of the semi-global, practical stability estimates that are obtained in extremum seeking problems and provides *checkable* conditions. Indeed, the calculation of the logarithmic norm of a square matrix is an easy task for various norms of $\mathbb{R}^n$. The existence of a function $p \in C^0(\mathbb{R}_+)$ and a norm $|\cdot|$ of $\mathbb{R}^n$ for which $\mu\left(\frac{\partial f}{\partial z}(t,z)\right) \leq p(t)$ for all $t \geq 0$, $z \in C$, where $\mu(\cdot)$ denotes the logarithmic norm associated with the norm $|\cdot|$, is guaranteed when there exists a symmetric, positive definite matrix $P \in \mathbb{R}^{n \times n}$ such that

$$\left(\frac{\partial f}{\partial z}(t,z)\right)^T P + P\left(\frac{\partial f}{\partial z}(t,z)\right) \leq 2p(t)P \text{ for all } t \geq 0, z \in C, \tag{2.9}$$

where $^T$ denotes transposition. Indeed, when $P \in \mathbb{R}^{n \times n}$ is a positive definite, symmetric matrix then there exists an invertible matrix $L \in \mathbb{R}^{n \times n}$ such that $P = L^T L$. Then the norm $|x| = |Lx|_2$, where $|\cdot|_2$ denotes the Euclidean norm in $\mathbb{R}^n$, gives that the logarithmic norm associated with the norm $|\cdot|$ of a matrix $G \in \mathbb{R}^{n \times n}$ is equal to the maximum eigenvalue of the symmetric matrix $\frac{1}{2}LGL^{-1} + \frac{1}{2}\left(L^{-1}\right)^T G^T L^T$. Consequently, $\max\left\{x^T\left(LGL^{-1} + \left(L^{-1}\right)^T G^T L^T\right)x : |x|_2 = 1\right\} = 2\mu(G)$, which implies $\max\left\{y^T\left(PG + G^T P\right)y : y^T P y = 1\right\} = 2\mu(G)$. Thus, the condition (2.9) implies the condition $\mu\left(\frac{\partial f}{\partial z}(t,z)\right) \leq p(t)$ for all $t \geq 0$, $z \in C$, where $\mu(\cdot)$ denotes the logarithmic norm associated with the norm $|\cdot|$.

## 3. Applications to Time-Invariant Systems

Theorem 1 has some interesting consequences to time-invariant (autonomous) systems. More specifically, we can obtain the result.

**Theorem 3:** *Consider the system*

$$\dot{x} = f(x), x \in \mathbb{R}^n \tag{3.1}$$

$$\dot{y} = g(x,y), y \in \mathbb{R} \tag{3.2}$$

*where $f \in C^1(\mathbb{R}^n; \mathbb{R}^n)$, $g \in C^1(\mathbb{R}^n \times \mathbb{R})$. Let $(x_0, y_0) \in \mathbb{R}^n \times \mathbb{R}$ be a point for which the solution $(x(t), y(t))$ of the initial value problem (3.1), (3.2) with initial condition $(x(0), y(0)) = (x_0, y_0)$ is bounded. Suppose that $\omega(x_0)$ contains a $T-$periodic orbit of system (3.1), where $\omega(x_0)$ is the omega limit set of $x_0 \in \mathbb{R}^n$ for system (3.1). Then $\omega(x_0, y_0)$ contains a $T-$periodic orbit of system (3.1), (3.2), where $\omega(x_0, y_0)$ is the omega limit set of $(x_0, y_0) \in \mathbb{R}^n \times \mathbb{R}$ for system (3.1), (3.2).*

The proof of Theorem 3 is provided in Section 5. Applying Theorem 3 inductively, we obtain the following result which deals with a cascade of time-invariant ODEs.

**Corollary 1:** *Consider system (3.1) with*

$$\dot{y}_i = g_i(x, y_1, ..., y_i), \; y_i \in \mathbb{R} \; \text{for} \; i = 1, ..., m \quad (3.3)$$

*where $f \in C^1(\mathbb{R}^n; \mathbb{R}^n)$, $g_i \in C^1(\mathbb{R}^n \times \mathbb{R}^i)$ for $i = 1, ..., m$. Let $(x_0, y_0) \in \mathbb{R}^n \times \mathbb{R}^m$ be a point for which the solution $(x(t), y_1(t), ..., y_m(t))$ of the initial value problem (3.1), (3.3) with initial condition $(x(0), y_1(0), ..., y_m(0)) = (x_0, y_0)$ is bounded. Suppose that $\omega(x_0)$ contains a $T$-periodic orbit of system (3.1), where $\omega(x_0)$ is the omega limit set of $x_0 \in \mathbb{R}^n$ for system (3.1). Then $\omega(x_0, y_0)$ contains a $T$-periodic orbit of system (3.1), (3.3), where $\omega(x_0, y_0)$ is the omega limit set of $(x_0, y_0) \in \mathbb{R}^n \times \mathbb{R}^m$ for system (3.1), (3.3).*

Proposition 2.3.3 on page 40 in [1] guarantees that the omega limit sets for system (3.1) with $n = 1$ can only contain equilibrium points. Moreover, the Poincaré-Bendixson theorem guarantees that the omega limit sets for system (3.1) with $n = 2$ are either empty or contain periodic orbits (notice that equilibria are periodic orbits). Therefore, we obtain from Corollary 1 the following corollary.

**Corollary 2:** *Consider system (3.1), (3.3) with $n = 1$ or $n = 2$, where $f \in C^1(\mathbb{R}^n; \mathbb{R}^n)$, $g_i \in C^1(\mathbb{R}^n \times \mathbb{R}^i)$ for $i = 1, ..., m$. Let $(x_0, y_0) \in \mathbb{R}^n \times \mathbb{R}^m$ be a point for which the solution $(x(t), y_1(t), ..., y_m(t))$ of the initial value problem (3.1), (3.3) with initial condition $(x(0), y_1(0), ..., y_m(0)) = (x_0, y_0)$ is bounded. Then $\omega(x_0, y_0)$ contains a periodic orbit of system (3.1), (3.3), where $\omega(x_0, y_0)$ is the omega limit set of $(x_0, y_0) \in \mathbb{R}^n \times \mathbb{R}^m$ for system (3.1), (3.3).*

It should be emphasized that neither Corollary 1 nor Corollary 2 guarantee that the omega limit set is a periodic orbit. The omega limit set may be a more complicated set (e.g. a homoclinic orbit).

**Example:** Consider the nonlinear system

$$\dot{x}_1 = x_2$$
$$\dot{x}_2 = \mu(1 - x_1^2)x_2 - x_1 \quad (3.4)$$
$$\dot{y} = -x_2^2 y + x_2^3$$

where $\mu > 0$ is a constant. The reader can see that the $x$-dynamics of system (3.4), i.e.,

$$\dot{x}_1 = x_2$$
$$\dot{x}_2 = \mu(1 - x_1^2)x_2 - x_1 \quad (3.5)$$

is the well-known van der Pol oscillator. Therefore, we know that for every $x_0 \in \mathbb{R}^2 \setminus \{0\}$ the solution of (3.5) with initial condition $x(0) = x_0$ is bounded and converges to a limit cycle of period $T > 0$. In other words, for every $x_0 \in \mathbb{R}^2 \setminus \{0\}$ $\omega(x_0)$ contains a $T$-periodic orbit of system (3.5), where $\omega(x_0)$ is the omega limit set of $x_0 \in \mathbb{R}^n$ for system (3.5).

The solutions of the ODE $\dot{y} = -x_2^2 y + x_2^3$ are bounded for every bounded and continuous function $x_2$. This follows from the fact that the following implication holds:

$$|y| \geq |x_2| \Rightarrow \frac{d}{dt}\left(\frac{1}{2} y^2\right) \leq 0$$

Therefore, for every $x_0 \in \mathbb{R}^2 \setminus \{0\}$, $y_0 \in \mathbb{R}$ the solution of (3.4) with initial condition $(x(0), y(0)) = (x_0, y_0)$ is bounded. Thus we can apply Theorem 3 or Corollary 1 and conclude that $\omega(x_0, y_0)$ contains a $T$-periodic orbit of system (3.4), where $\omega(x_0, y_0)$ is the omega limit set of $(x_0, y_0) \in \mathbb{R}^2 \times \mathbb{R}$ for system (3.4). A straightforward contradiction argument can guarantee that the $T$-periodic orbit of system (3.4) contained in $\omega(x_0, y_0)$ is not an equilibrium point.

To the best of our knowledge such a result cannot be guaranteed by existing results in the literature. ◁

## 4. Applications to Extremum-Seeking

Theorem 1 and Theorem 2 have some interesting consequences to extremum-seeking problems. Consider the following system

$$\dot{x}(t) = -\varepsilon \sin(t) h(x(t) + a \sin(t)) \quad (4.1)$$

where $\varepsilon, a > 0$ are constants, $x(t) \in \mathbb{R}$ is the state and $h \in C^3(\mathbb{R})$ is a function that satisfies the following properties:

$$h'(x) x > 0, \text{ for all } x \neq 0 \quad (4.2)$$

$$h'(0) = 0, \quad h''(0) > 0 \quad (4.3)$$

Clearly, system (4.1) is a system of the form (2.1) with $n = 1$, $f(t, x) = -\varepsilon \sin(t) h(x + a \sin(t))$ (even without assuming (4.2), (4.3)) and Theorem 1 can be applied to system (4.1) with $T = 2\pi$: for every bounded solution $x(t) \in \mathbb{R}$ of (4.1) there exists a $2\pi$-periodic solution $x^*(t)$ of (4.1) such that $\lim_{t \to +\infty}(|x(t) - x^*(t)|) = 0$.

System (4.1) under assumptions (4.2), (4.3) arises in extremum seeking for the static output map $y = h(x)$ with no plant dynamics; see for instance [10, 11, 12, 14]. Here, the extremum is taken at $x = 0$ for notational simplicity and without loss of generality. The reader can immediately understand the novel perspective that Theorem 1 offers for system (4.1): here we do not need to apply averaging arguments (as in [6, 7]) and the results are not local.

System (4.1) under assumptions (4.2), (4.3) was also studied in [18]; see Corollary 1 in [18] where it is shown that system (4.1) is Semi-globally Practically Asymptotically Stable, i.e., there exists a function $\beta \in KL$ for which the following property holds:

**(P)** *For every* $R, \kappa > 0$, *there exist constants* $\varepsilon^*, a^* > 0$ *(depending on* $R, \kappa > 0$*) such that for every fixed* $\varepsilon \in (0, \varepsilon^*)$, $a \in (0, a^*)$ *and for every* $t_0 \geq 0$, $x_0 \in \mathbb{R}$ *with* $|x_0| \leq R$, *the unique solution* $x(t) \in \mathbb{R}$ *of (4.1) with initial condition* $x(t_0) = x_0$ *exists for all* $t \geq t_0$ *and satisfies the following estimate for all* $t \geq t_0$:

$$|x(t)| \leq \beta\big(|x_0|, \varepsilon a(t-t_0)\big) + \kappa \tag{4.4}$$

Of course, property (P) can be combined with Theorem 1 but that's not the end of the story. We want to exploit Theorem 2 which can be applied if there exist a constant $b > 0$ and a $2\pi$-periodic function $p \in C^0(\mathbb{R}_+)$ with $\int_0^{2\pi} p(s)ds < 0$ such that $\frac{\partial f}{\partial z}(t,z) = -\varepsilon \sin(t) h'(z + a\sin(t)) \leq p(t)$ for all $t \geq 0$, $z \in \mathbb{R}$ with $|z| \leq b$ (that is we use the set $C = [-b,b]$ in Theorem 2). Clearly, if

$$4b + 2\frac{(b+a)^2}{h''(0)} \max_{|s| \leq b+a}\big(|h'''(s)|\big) < a\pi \tag{4.5}$$

then the aforementioned assumption of Theorem 2 is valid with

$$p(t) := -\varepsilon a h''(0) \sin^2(t) + \varepsilon b |\sin(t)| h''(0) + \frac{\varepsilon}{2} |\sin(t)| (b+a)^2 \max_{|s| \leq b+a}\big(|h'''(s)|\big) \tag{4.6}$$

Definition (4.6), inequalities (4.3), (4.5) and the fact that $\int_0^{2\pi} |\sin(s)| ds = 2\int_0^{\pi} \sin(s) ds = 4$ guarantee that $\int_0^{2\pi} p(s) ds < 0$ and we also have for all $t \geq 0$, $z \in \mathbb{R}$ with $|z| \leq b$:

$$\frac{\partial f}{\partial z}(t,z) = -\varepsilon \sin(t) h'(z + a\sin(t)) = -\varepsilon \sin(t)\big(h'(z + a\sin(t)) - h'(0)\big)$$

$$= -\varepsilon \sin(t)(z + a\sin(t))\int_0^1 h''(\lambda z + a\lambda \sin(t)) d\lambda$$

$$= -\varepsilon \sin(t)(z + a\sin(t)) h''(0)$$

$$-\varepsilon \sin(t)(z + a\sin(t))\int_0^1 \big(h''(\lambda z + a\lambda \sin(t)) - h''(0)\big) d\lambda$$

$$= -\varepsilon \sin(t) h''(0) z - \varepsilon a h''(0) \sin^2(t)$$

$$-\varepsilon \sin(t)(z + a\sin(t))^2 \int_0^1\int_0^1 \lambda h'''(\lambda \mu z + a\lambda\mu \sin(t)) d\mu d\lambda$$

$$\leq -\varepsilon a h''(0) \sin^2(t) + \varepsilon b |\sin(t)| h''(0) + \frac{\varepsilon}{2} |\sin(t)|(b+a)^2 \max_{|s| \leq b+a}\big(|h'''(s)|\big) = p(t)$$

Consequently, in order to be able to apply Theorem 2 to system (4.1) we need to find a set $S \subseteq \mathbb{R}$ and a constant $b > 0$ for which (4.5) holds such that $\limsup_{t \to +\infty}(|x(t)|) < b$ for every solution $x(t)$ of (4.1) with $x(t_0) \in S$ for some $t_0 \geq 0$. To this purpose, we need a sharp estimate of the way that the residual constant $\kappa > 0$ in (4.4) depends on the parameters $\varepsilon, a > 0$. This sharp estimate is provided by the lemma below. The proof of the following lemma is based on the time-periodic transformation which is used in the proof of the averaging theorem (Theorem 8.3 in [5]) and was also used in the paper [18].

**Lemma 1:** *Consider system (4.1) under assumptions (4.2), (4.3) with $\varepsilon > 0$, $a \in (0,1]$. Let arbitrary $R > 0$ be given and suppose that the following inequalities hold:*

$$8\varepsilon \max_{|z| \leq 2R+2}(|h'(z)|) < 1 \tag{4.7}$$

$$8\varepsilon \max_{|z| \leq 2R+1}(|h(z)|) \leq R \tag{4.8}$$

$$\frac{a^2}{6\sigma} \max_{|z| \leq 2R+2}(|h'''(z)|) + \frac{16\varepsilon}{a\sigma} \max_{|z| \leq 2R+2}(|h(z)|) \max_{|z| \leq 3R+2}(|h'(z)|) \leq R \tag{4.9}$$

*where $\sigma > 0$ is defined by the following equation:*

$$\sigma := \inf\left\{\frac{wh'(w)}{w^2} : 0 < |w| \leq 2R+1\right\} = \min\left\{\int_0^1 h''(\lambda w) d\lambda : |w| \leq 2R+1\right\} \tag{4.10}$$

*Then for every $t_0 \geq 0$, $x_0 \in \mathbb{R}$ with $|x_0| \leq R$ the unique solution $x(t)$ of (4.1) with initial condition $x(t_0) = x_0$ is defined for $t \geq t_0$ and satisfies*

$$\limsup_{t \to +\infty}(|x(t)|) \leq \frac{a^2}{6\sigma} \max_{|z| \leq 2R+2}(|h'''(z)|) + 8\varepsilon \max_{|z| \leq 2R+2}(|h(z)|)\left(\frac{2}{a\sigma} \max_{|z| \leq 3R+2}(|h'(z)|) + 1\right) \tag{4.11}$$

**Remarks:** (i) The fact that $\sigma$ defined by (4.10) is positive is a consequence of (4.2) and (4.3).
(ii) There is a clear qualitative similarity between (4.11) and estimate (37) in [7]. However, there is a difference between (4.11) and estimate (37) in [7]: while estimate (37) in [7] is local, estimate (4.11) is not local. Estimate (4.11) is a regional estimate, in the sense that it holds for all initial conditions in the region (interval) $[-R, R]$. From a different point of view estimate (4.11) is semi-global, because for every $R > 0$ inequalities (4.7), (4.8), (4.9) (and consequently (4.11) for all initial conditions in $[-R, R]$) can be guaranteed for sufficiently small $\varepsilon > 0$ and $a \in (0,1]$.

We are now ready to combine Theorem 1, Theorem 2 and Lemma 1.

**Theorem 4:** *Consider system (4.1) under assumptions (4.2), (4.3) with $\varepsilon, a > 0$. Then for every bounded solution $x(t)$ of (4.1) there exists a $2\pi$-periodic solution $x^*(t)$ of (4.1) such that $\lim_{t \to +\infty} (|x(t) - x^*(t)|) = 0$. Furthermore, if $a \leq 1$ and if the following inequalities hold for certain $R > 0$*

$$a\gamma \left( 12\sigma h''(0) + (\gamma + 6\sigma)^2 \right) < 18\sigma^2 h''(0) \pi \tag{4.12}$$

$$a^2 \gamma < 6\sigma R \tag{4.13}$$

*where $\sigma > 0$ is defined by (4.10) and $\gamma := \max_{|z| \leq 2R+2} (|h'''(z)|)$, then there exists $\varepsilon^* > 0$ such that the following property holds for every $\varepsilon \in (0, \varepsilon^*)$:*

**(P')** *There exists a unique $2\pi$-periodic solution $x^{**}(t)$ of (4.1) with $\max_{t \in [0, 2\pi]} (|x^{**}(t)|) \leq \dfrac{a^2 \gamma}{6\sigma} + \varepsilon K$, where $K := 8 \max_{|z| \leq 2R+2} (|h(z)|) \left( \dfrac{2}{a\sigma} \max_{|z| \leq 3R+2} (|h'(z)|) + 1 \right)$, such that for every $t_0 \geq 0$, $x_0 \in \mathbb{R}$ with $|x_0| \leq R$ the unique solution $x(t)$ of (4.1) with initial condition $x(t_0) = x_0$ is defined for $t \geq t_0$ and satisfies $\lim_{t \to +\infty} (|x(t) - x^{**}(t)|) = 0$.*

The reader can easily spot the differences with extremum-seeking results that make use of the averaging theorem: the extremum-seeking results that make use of the averaging theorem are local. On the other hand, Theorem 4 is inherently non-local. Theorem 4 provides a complete description of the dynamics of (4.1): either a solution of (4.1) is bounded and converges to a periodic solution or is unbounded. Moreover, there exists a unique periodic solution $x^{**}(t)$ of (4.1) which is located close to zero (this is guaranteed by the estimate $\max_{t \in [0, 2\pi]} (|x^{**}(t)|) \leq \dfrac{a^2 \gamma}{6\sigma} + \varepsilon K$ and this feature is also guaranteed in [7]) that attracts the solutions initiated in the interval $[-R, R]$ (a feature absent from the results in [7]).

For the quadratic case, i.e., when $h(x) = \vartheta + Ex^2$ with $E > 0$ we get $\sigma = 2E$ and $\gamma = 0$. It follows that in the quadratic case inequalities (4.12), (4.13) are automatically valid for every $R > 0$. We conclude from Theorem 4 that for every $a \leq 1$, $R > 0$ there exists $\varepsilon^* > 0$ such that the following property holds for every $\varepsilon \in (0, \varepsilon^*)$:

**(P'')** *There exists a unique $2\pi$-periodic solution $x^{**}(t)$ of (4.1) with $h(x) = \vartheta + Ex^2$ satisfying $\max_{t \in [0, 2\pi]} (|x^{**}(t)|) \leq 8\varepsilon \dfrac{6R + 4 + a}{a} \left( |\vartheta| + 4E(R+1)^2 \right)$, such that for every $t_0 \geq 0$, $x_0 \in \mathbb{R}$ with $|x_0| \leq R$ the unique solution $x(t)$ of (4.1) with $h(x) = \vartheta + Ex^2$ and initial condition $x(t_0) = x_0$ is defined for $t \geq t_0$ and satisfies $\lim_{t \to +\infty} (|x(t) - x^{**}(t)|) = 0$.*

Moreover, in the quadratic case *there is no bias*, in the sense that the $2\pi$-periodic solution $x^{**}(t)$ passes through zero (at least) once in every interval of length $2\pi$ and therefore, *zero is a limit point* of all solutions initiated in the interval $[-R, R]$. This is a result of the following proposition which is a direct consequence of Banach's fixed point theorem and the symmetry introduced by an even function $h$.

**Proposition 1:** *Suppose that $h \in C^3(\mathbb{R})$ is an even function. Let $R > 0$ be given and suppose that*

$$\frac{\varepsilon\pi}{2} \max_{|z| \leq R+a} (|h(z)|) \leq R \text{ and } \frac{\varepsilon\pi}{2} \max_{|z| \leq R+a} (|h'(z)|) < 1 \tag{4.14}$$

*Then there exists a $2\pi$-periodic solution $x^{**}$ of (4.1) with $x^{**}(0) = -x^{**}(\pi)$ and $\max_{t \in [0, 2\pi]} (|x^{**}(t)|) \leq \frac{\varepsilon\pi}{2} \max_{|z| \leq R+a} (|h(z)|)$.*

Proposition 1 provides a sharp estimate of the magnitude of the $2\pi$-periodic solution $x^{**}$ of (4.1) that attracts all solutions initiated in the interval $[-R, R]$.

## 5. Proofs

In this section we provide all proofs of the results that were stated above. We start with the proof of Theorem 1.

**Proof of Theorem 1:** We define the following set for each $\tau \geq 0$:

$$D(\tau) := \left\{ z \in \mathbb{R} : \sup_{t \geq \tau} (|\varphi(t, \tau, z)|) < +\infty \right\} \tag{5.1}$$

Definition (5.1) in conjunction with Fact II implies that $D(\tau) \subseteq \mathbb{R}$ is an interval (a connected set) with $\varphi(t, \tau, D(\tau)) \subseteq D(t)$ for all $t \geq \tau \geq 0$ (a consequence of the semigroup property). Indeed, if $x, y \in D(\tau)$ with $x \leq y$ then for every $z \in [x, y]$ we have (due to Fact II) $\varphi(t, \tau, z) \in [\varphi(t, \tau, x), \varphi(t, \tau, y)]$ for all $t \geq \tau$ for which $\varphi(t, \tau, z)$ exists. Since both $\varphi(t, \tau, x), \varphi(t, \tau, y)$ are bounded we obtain that $\varphi(t, \tau, z)$ exists for all $t \geq \tau$ and is bounded. Fact I implies that $D(\tau + kT) = D(\tau)$, for all $k \in \mathbb{N}$, $\tau \geq 0$. Furthermore, definition (5.1) and the semigroup property imply that the existence of $t \geq \tau$ with $\varphi(t, \tau, z) \in D(t)$ guarantees that $z \in D(\tau)$.

We also define the set for each $\tau, r \geq 0$:

$$D_r(\tau) := \left\{ z \in \mathbb{R} : \sup_{t \geq \tau}\left(|\varphi(t,\tau,z)|\right) \leq r \right\} \quad (5.2)$$

When the set $D_r(\tau) \subseteq \mathbb{R}$ is non-empty, it is an interval (a connected set) with $\varphi(t,\tau,D_r(\tau)) \subseteq D_r(t)$ for all $t \geq \tau \geq 0$ (a consequence of the semigroup property). Indeed, if $x, y \in D_r(\tau)$ with $x \leq y$ then for every $z \in [x,y]$ we have (due to Fact II) $\varphi(t,\tau,z) \in [\varphi(t,\tau,x), \varphi(t,\tau,y)] \subseteq [-r,r]$ for all $t \geq \tau$ for which $\varphi(t,\tau,z)$ exists. Therefore $\varphi(t,\tau,z)$ exists for all $t \geq \tau$ and $\varphi(t,\tau,z) \in [-r,r]$ for all $t \geq \tau$. Consequently, definition (5.2) implies that $z \in D_r(\tau)$.

The set $D_r(\tau) \subseteq \mathbb{R}$ is closed. Clearly, this is valid when $D_r(\tau) = \varnothing$. When $D_r(\tau) \neq \varnothing$ then we can consider a sequence $\{z_k \in D_r(\tau) : k \in \mathbb{N}\}$ with $\lim(z_k) = z^* \in \mathbb{R}$. We show next that $z^* \in D_r(\tau)$. Clearly, there exists $t_{\max} \in (\tau, +\infty]$ for which the solution $\varphi(t,\tau,z^*)$ is defined for $t \in [\tau, t_{\max})$. Moreover, if $t_{\max} < +\infty$ then $\limsup_{t \to t_{\max}^-}\left(|\varphi(t,\tau,z^*)|\right) = +\infty$. For every $t \in [\tau, t_{\max})$ we get that $|\varphi(t,\tau,z^*)| \leq |\varphi(t,\tau,z_k)| + |\varphi(t,\tau,z^*) - \varphi(t,\tau,z_k)|$ for all $k \in \mathbb{N}$. Since $z_k \in D_r(\tau)$ for all $k \in \mathbb{N}$, it follows from definition (5.2) that $|\varphi(t,\tau,z_k)| \leq r$ for all $k \in \mathbb{N}$. Therefore, $|\varphi(t,\tau,z^*)| \leq r + |\varphi(t,\tau,z^*) - \varphi(t,\tau,z_k)|$ for all $k \in \mathbb{N}$. Since $\lim(z_k) = z^*$, by continuity of the solution map we get $\lim(\varphi(t,\tau,z_k)) = \varphi(t,\tau,z^*)$, which implies that $|\varphi(t,\tau,z^*)| \leq r$ for all $t \in [\tau, t_{\max})$. This implies that $t_{\max} = +\infty$ and $\varphi(t,\tau,z^*) \in [-r,r]$ for $t \geq \tau$. Consequently, definition (5.2) implies that $z^* \in D_r(\tau)$.

Define the map $g : D(t_0) \to D(t_0)$ by means of the following formula for $x \in D(t_0)$:

$$g(x) = \varphi(t_0 + T, t_0, x) \quad (5.3)$$

The map $g : D(t_0) \to D(t_0)$ defined by (5.3) is the so-called first return (Poincaré) map for (2.1). The map $g : D(t_0) \to D(t_0)$ is continuous (by continuity of the solution map) and due to Fact II, is increasing.

Since $x_0 \in D(t_0)$ we define $r := \sup_{t \geq t_0}\left(|\varphi(t,t_0,x_0)|\right)$. Therefore, definition (5.2) implies that $x_0 \in D_r(t_0)$. We define for all $k \in \mathbb{N}$:

$$y_k = \varphi(t_0 + kT, t_0, x_0) \quad (5.4)$$

Notice that definition (5.4) and the fact that $\varphi(t,t_0,D_r(t_0)) \subseteq D_r(t)$ for all $t \geq t_0 \geq 0$ implies that $y_k \in D_r(t_0)$ for all $k \in \mathbb{N}$. Moreover, using the semigroup property, definition (5.3) and Fact I, we get:

$$y_{k+1} = \varphi(t_0 + (k+1)T, t_0, x_0)$$
$$= \varphi(t_0 + (k+1)T, t_0 + kT, \varphi(t_0 + kT, t_0, x_0)) \quad (5.5)$$
$$= \varphi(t_0 + (k+1)T, t_0 + kT, y_k) = \varphi(t_0 + T, t_0, y_k) = g(y_k)$$

If $g(y_0) = y_0$ then $y_k \equiv y_0$. Since $y_0 = x_0$, we get from definition (5.3) that $\varphi(t_0 + T, t_0, x_0) = x_0$, which implies that the solution $\varphi(t, t_0, x_0)$ is periodic of period $T > 0$.

If $y_1 = g(y_0) > y_0$ then the fact that $g$ is increasing implies that the sequence $\{y_k \in D_r(t_0) : k \in \mathbb{N}\}$ is increasing with $|y_k| \le r$ for all $k \in \mathbb{N}$. Therefore, there exists $y^* \in \mathbb{R}$ such that $y^* = \lim(y_k)$. Since the set $D_r(t_0) \subseteq \mathbb{R}$ is closed, it follows that $y^* \in D_r(t_0) \subseteq D(t_0)$. Therefore, continuity of the map $g : D(t_0) \to D(t_0)$ implies that $y^* = \lim(y_{k+1}) = \lim(g(y_k)) = g(\lim(y_k)) = g(y^*)$. Definition (5.3) that $\varphi(t_0 + T, t_0, y^*) = y^*$, which implies that the solution $\varphi(t, t_0, y^*)$ is periodic of period $T > 0$.

Similarly, if $y_1 = g(y_0) < y_0$ then the fact that $g$ is increasing implies that the sequence $\{y_k \in D_r(t_0) : k \in \mathbb{N}\}$ is decreasing with $|y_k| \le r$ for all $k \in \mathbb{N}$. Therefore, there exists $y^* \in \mathbb{R}$ such that $y^* = \lim(y_k)$. Since the set $D_r(t_0) \subseteq \mathbb{R}$ is closed, it follows that $y^* \in D_r(t_0) \subseteq D(t_0)$. Therefore, continuity of the map $g : D(t_0) \to D(t_0)$ implies that $y^* = \lim(y_{k+1}) = \lim(g(y_k)) = g(\lim(y_k)) = g(y^*)$. Definition (5.3) that $\varphi(t_0 + T, t_0, y^*) = y^*$, which implies that the solution $\varphi(t, t_0, y^*)$ is periodic of period $T > 0$.

We next show that $\lim_{t \to +\infty} \left( \left| \varphi(t, t_0, x_0) - \varphi(t, t_0, y^*) \right| \right) = 0$. Let $k \in \mathbb{N}$ and $t \in [t_0 + kT, t_0 + (k+1)T)$ be given. Clearly, we get:

$$k = \left\lfloor \frac{t - t_0}{T} \right\rfloor \quad (5.6)$$

Definition (5.6) implies that
$$\tau \in [0, T) \quad (5.7)$$
where
$$\tau = t - t_0 - kT \quad (5.8)$$

Using the semigroup property, definitions (5.4), (5.8) and Fact I, we get:

$$\varphi(t, t_0, x_0) = \varphi(\tau + t_0 + kT, t_0, x_0)$$
$$= \varphi(\tau + t_0 + kT, t_0 + kT, \varphi(t_0 + kT, t_0, x_0)) \quad (5.9)$$
$$= \varphi(\tau + t_0 + kT, t_0 + kT, y_k) = \varphi(\tau + t_0, t_0, y_k)$$

Using the semigroup property, definition (5.8), Fact I and the fact that $\varphi(t_0 + T, t_0, y^*) = y^*$, we get:

$$\varphi(t,t_0,y^*) = \varphi(\tau+t_0+kT,t_0,y^*)$$
$$= \varphi(\tau+t_0+kT,t_0+kT,\varphi(t_0+kT,t_0,y^*)) \qquad (5.10)$$
$$= \varphi(\tau+t_0+kT,t_0+kT,y^*) = \varphi(\tau+t_0,t_0,y^*)$$

Combining (5.9), (5.10), (2.3) and the triangle inequality we obtain:

$$\left|\varphi(t,t_0,x_0) - \varphi(t,t_0,y^*)\right| = \left|\varphi(\tau+t_0,t_0,y_k) - \varphi(\tau+t_0,t_0,y^*)\right|$$
$$= \left|(y_k - y^*) + \int_{t_0}^{\tau+t_0}\left(f(s,\varphi(s,t_0,y_k)) - f(s,\varphi(s,t_0,y^*))\right)ds\right| \qquad (5.11)$$
$$\leq \left|y_k - y^*\right| + \int_{t_0}^{\tau+t_0}\left|f(s,\varphi(s,t_0,y_k)) - f(s,\varphi(s,t_0,y^*))\right|ds$$

Since $y^* \in D_r(t_0)$ and $y_k \in D_r(t_0)$, definition (5.2) implies that the following estimate holds for all $s \in [t_0, t_0+T]$:

$$\left|f(s,\varphi(s,t_0,y_k)) - f(s,\varphi(s,t_0,y^*))\right| \leq L\left|\varphi(s,t_0,y_k) - \varphi(s,t_0,y^*)\right| \qquad (5.12)$$

where $L = \max\left\{\left|\dfrac{\partial f}{\partial x}(s,z)\right| : s \in [t_0,t_0+T], |z| \leq r\right\}$. Therefore, we obtain from (5.7), (5.11) and (5.12) for all $\tau \in [0,T]$:

$$\left|\varphi(\tau+t_0,t_0,y_k) - \varphi(\tau+t_0,t_0,y^*)\right|$$
$$\leq \left|y_k - y^*\right| + L\int_{t_0}^{\tau+t_0}\left|\varphi(s,t_0,y_k) - \varphi(s,t_0,y^*)\right|ds \qquad (5.13)$$
$$= \left|y_k - y^*\right| + L\int_{0}^{\tau}\left|\varphi(s+t_0,t_0,y_k) - \varphi(s+t_0,t_0,y^*)\right|ds$$

Using (5.13) and the Gronwall-Bellman inequality, we obtain for all $\tau \in [0,T]$:

$$\left|\varphi(\tau+t_0,t_0,y_k) - \varphi(\tau+t_0,t_0,y^*)\right|$$
$$\leq \exp(L\tau)\left|y_k - y^*\right| \leq \exp(LT)\left|y_k - y^*\right| \qquad (5.14)$$

Combining (5.14), (5.6) and (5.11) we get for all $t \geq t_0$:

$$\left|\varphi(t,t_0,x_0) - \varphi(t,t_0,y^*)\right| \leq \exp(LT)\left|y_{\lfloor(t-t_0)/T\rfloor} - y^*\right| \qquad (5.15)$$

Estimate (5.15) and the fact that $y^* = \lim(y_k)$ imply that $\lim_{t \to +\infty} \left( \left| \varphi(t, t_0, x_0) - \varphi(t, t_0, y^*) \right| \right) = 0$. The proof is complete. ◁

For the proof of Theorem 2 we need the following technical lemma.

**Lemma 2:** *Let $A : [0,1] \to \mathbb{R}^{n \times n}$ be a continuous map. Let $|\cdot|$ be a norm of $\mathbb{R}^n$ and let $\mu(\cdot)$ denote the logarithmic norm of a $n \times n$ matrix associated with the norm $|\cdot|$. Suppose that there exists a constant $p \in \mathbb{R}$ such that $\mu(A(\lambda)) \leq p$ for all $\lambda \in [0,1]$. Then $\mu\left( \int_0^1 A(\lambda) d\lambda \right) \leq p$.*

**Proof:** Let $\|\cdot\|$ be the matrix norm of $\mathbb{R}^{n \times n}$ induced by the norm $|\cdot|$. Clearly, there exists a constant $K > 0$ such that for every matrix $B = \{b_{i,j} : i, j = 1,...,n\} \in \mathbb{R}^{n \times n}$ it holds that

$$\|B\| \leq K \max_{i,j=1,...,n} \left( |b_{i,j}| \right) \tag{5.16}$$

Set $A(\lambda) = \{a_{i,j}(\lambda) : i, j = 1,...,n\} \in \mathbb{R}^{n \times n}$ where $a_{i,j} : [0,1] \to \mathbb{R}$ are continuous functions for $i, j = 1,...,n$.

Let arbitrary $\varepsilon > 0$ be given. Since $a_{i,j} : [0,1] \to \mathbb{R}$ are continuous functions for $i, j = 1,...,n$ there exist numbers $0 = t_0 < t_1 < ... < t_m = 1$ such that $\left| \int_0^1 a_{i,j}(\lambda) d\lambda - \sum_{k=1}^m (t_k - t_{k-1}) a_{i,j}(t_k) \right| < \varepsilon$ for all $i, j = 1,...,n$. We define the matrix $G := \left\{ g_{i,j} = \int_0^1 a_{i,j}(\lambda) d\lambda - \sum_{k=1}^m (t_k - t_{k-1}) a_{i,j}(t_k) : i, j = 1,...,n \right\}$ and notice that due to (5.16) we obtain:

$$\int_0^1 A(\lambda) d\lambda = G + \sum_{k=1}^m (t_k - t_{k-1}) A(t_k) \tag{5.17}$$

$$\|G\| \leq K\varepsilon \tag{5.18}$$

Exploiting the fact that $\mu(C + D) \leq \mu(C) + \mu(D)$, $\mu(C) \leq \|C\|$ and $\mu(\gamma C) = \gamma \mu(C)$ for every pair of matrices $C, D \in \mathbb{R}^{n \times n}$ and for every $\gamma > 0$, we obtain from (5.17), (5.18) and the fact that $\mu(A(\lambda)) \leq p$ for all $\lambda \in [0,1]$:

$$\mu\left(\int_0^1 A(\lambda)d\lambda\right) = \mu\left(G + \sum_{k=1}^m (t_k - t_{k-1})A(t_k)\right) \leq \mu(G) + \sum_{k=1}^m \mu\left((t_k - t_{k-1})A(t_k)\right)$$

$$\leq \|G\| + \sum_{k=1}^m (t_k - t_{k-1})\mu(A(t_k)) \leq K\varepsilon + p\sum_{k=1}^m (t_k - t_{k-1}) = K\varepsilon + p$$

Since the above inequality holds for arbitrary $\varepsilon > 0$, we obtain that $\mu\left(\int_0^1 A(\lambda)d\lambda\right) \leq p$. The proof is complete. ◁

We continue with the proof of Theorem 2 and the proof of Theorem 3.

**Proof of Theorem 2:** Let arbitrary $\tau_0 \geq 0$, $\xi \in S$ be given. Since there exists $\tau \geq \tau_0$ with $\varphi(t, \tau_0, \xi) \in C$ for all $t \geq \tau$, it follows from Fact I that the set

$$\bar{D}(l) := \{z \in \mathbb{R}^n : \varphi(t, l, z) \in C, \forall t \geq l\} \tag{5.19}$$

is non-empty for every $l \geq 0$. Moreover, it holds that $\varphi(t, l, \bar{D}(l)) \subseteq \bar{D}(l)$ for all $t \geq l \geq 0$ (a consequence of the semigroup property). Fact I implies that $\bar{D}(l + kT) = \bar{D}(l)$, for all $k \in \mathbb{N}$, $l \geq 0$.

The set $\bar{D}(\tau) \subseteq \mathbb{R}^n$ is closed. Indeed, we can consider a sequence $\{z_k \in \bar{D}(\tau) : k \in \mathbb{N}\}$ with $\lim(z_k) = z^* \in \mathbb{R}^n$. We show next that $z^* \in \bar{D}(\tau)$. Clearly, there exists $t_{\max} \in (\tau, +\infty]$ for which the solution $\varphi(t, \tau, z^*)$ is defined for $t \in [\tau, t_{\max})$. Moreover, if $t_{\max} < +\infty$ then $\limsup_{t \to t_{\max}^-}(\|\varphi(t, \tau, z^*)\|) = +\infty$. Since $z_k \in \bar{D}(\tau)$ for all $k \in \mathbb{N}$, it follows from definition (5.19) that $\varphi(t, \tau, z_k) \in C$ for all $k \in \mathbb{N}$. Since $\lim(z_k) = z^*$, by continuity of the solution map we get $\lim(\varphi(t, \tau, z_k)) = \varphi(t, \tau, z^*)$ and closure of $C \subseteq \mathbb{R}^n$ implies that $\varphi(t, \tau, z^*) \in C$ for all $t \in [\tau, t_{\max})$. Boundedness of $C \subseteq \mathbb{R}^n$ implies that $t_{\max} = +\infty$ and consequently $\varphi(t, \tau, z^*) \in C$ for all $t \geq \tau$. Consequently, definition (5.19) implies that $z^* \in \bar{D}(\tau)$.

Let arbitrary $x_0, y_0 \in \bar{D}(t_0)$ be given. Using (2.1) we obtain for all $t \geq t_0$:

$$\frac{\partial}{\partial t}(\varphi(t, t_0, x_0) - \varphi(t, t_0, y_0)) = f(t, \varphi(t, t_0, x_0)) - f(t, \varphi(t, t_0, y_0))$$
$$= \left(\int_0^1 \frac{\partial f}{\partial x}(t, (1-\lambda)\varphi(t, t_0, y_0) + \lambda\varphi(t, t_0, x_0))d\lambda\right)(\varphi(t, t_0, x_0) - \varphi(t, t_0, y_0)) \tag{5.20}$$

Exploiting the Coppel inequality (see [19]) and (5.20) we obtain for all $t \geq t_0$:

$$|\varphi(t,t_0,x_0)-\varphi(t,t_0,y_0)|$$
$$\leq |x_0-y_0|\exp\left(\int_{t_0}^{t}\mu\left(\int_0^1 \frac{\partial f}{\partial x}(s,(1-\lambda)\varphi(s,t_0,y_0)+\lambda\varphi(s,t_0,x_0))d\lambda\right)ds\right) \quad (5.21)$$

Since $\varphi(t,t_0,x_0)\in C$ and $\varphi(t,t_0,y_0)\in C$ for all $t\geq t_0$ (recall definition (5.19) and the fact that $x_0,y_0 \in \bar{D}(t_0)$), due to convexity of the set $C\subseteq \mathbb{R}^n$ it holds that $((1-\lambda)\varphi(s,t_0,y_0)+\lambda\varphi(s,t_0,x_0))\in C$ for all $\lambda\in[0,1]$ and $s\geq t_0$. Since $\mu\left(\frac{\partial f}{\partial z}(t,z)\right)\leq p(t)$ for all $t\geq 0$, $z\in C$, it follows that $\mu\left(\frac{\partial f}{\partial x}(s,(1-\lambda)\varphi(s,t_0,y_0)+\lambda\varphi(s,t_0,x_0))\right)\leq p(s)$ for all $\lambda\in[0,1]$ and $s\geq t_0$. Exploiting Lemma 2 and (5.21) we get for all $t\geq t_0$:

$$|\varphi(t,t_0,x_0)-\varphi(t,t_0,y_0)|\leq |x_0-y_0|\exp\left(\int_{t_0}^{t}p(s)ds\right) \quad (5.22)$$

Thus
$$|\varphi(t_0+T,t_0,x_0)-\varphi(t_0+T,t_0,y_0)|\leq |x_0-y_0|\exp(-c) \quad (5.23)$$

where $c:=-\int_0^T p(s)ds > 0$. Define the map $g:\bar{D}(t_0)\to \bar{D}(t_0+T)=\bar{D}(t_0)$ by means of formula (5.3) for $x\in\bar{D}(t_0)$. The map $g:\bar{D}(t_0)\to\bar{D}(t_0)$ defined by (5.3) is the return (Poincaré) map for (2.1). Due to (5.23) the map $g:\bar{D}(t_0)\to\bar{D}(t_0)$ is a contraction. Banach's fixed point theorem implies the existence of a fixed point $y^*\in\bar{D}(t_0)$ with $\varphi(t_0+T,t_0,y^*)=y^*$, which implies that the solution $x^*(t)=\varphi(t,t_0,y^*)$ is periodic of period $T>0$. Due to (2.2) the periodic solution $x^*(t)$ can be defined for all $t\geq 0$ and satisfies $x^*(t)\in C$ for all $t\geq 0$.

We next show that $\lim_{t\to+\infty}\left(|\varphi(t,t_0,x_0)-x^*(t)|\right)=0$ for all $t_0\geq 0$ and $x_0\in S$.

Let arbitrary $t_0\geq 0$, $x_0\in S$ be given. It follows that there exists $\tau>t_0$ such that $\varphi(t,t_0,x_0)\in C$ for all $t\geq\tau$. Definition (5.19) and the semigroup property imply that $\varphi(t,t_0,x_0)\in\bar{D}(t)$ and $x^*(t)=\varphi(t,t_0,y^*)\in\bar{D}(t)$ for all $t\geq\tau$.

Working in the same way as above we obtain for all $t\geq\bar{\tau}\geq\tau$:

$$|\varphi(t,t_0,x_0)-x^*(t)|\leq |\varphi(\bar{\tau},t_0,x_0)-x^*(\bar{\tau})|\exp\left(\int_{\bar{\tau}}^{t}p(s)ds\right) \quad (5.24)$$

Let $m \in \mathbb{N}$ with $mT \geq \tau$ be given. Applying (5.24) with $\bar{\tau} = (m+k)T$, $t = (m+k+1)T$, where $k \in \mathbb{N}$, we obtain for all $k \in \mathbb{N}$:

$$\left|\varphi\big((m+k+1)T,t_0,x_0\big) - x^*(0)\right|$$
$$\leq \left|\varphi\big((m+k)T,t_0,x_0\big) - x^*(0)\right| \exp\left(\int_{(m+k)T}^{(m+k+1)T} p(s)\,ds\right) \tag{5.25}$$
$$= \left|\varphi\big((m+k)T,t_0,x_0\big) - x^*(0)\right| \exp(-c)$$

where $c := -\int_0^T p(s)\,ds > 0$. Applying induction and using (5.25) we get for all $k \in \mathbb{N}$:

$$\left|\varphi\big((m+k)T,t_0,x_0\big) - x^*(0)\right| \leq \left|\varphi(mT,t_0,x_0) - x^*(0)\right| \exp(-kc) \tag{5.26}$$

For every $t \geq mT$, applying (5.24) with $\bar{\tau} = (m+k)T$, $k = \left\lfloor \dfrac{t}{T} \right\rfloor - m$, exploiting periodicity of $p$ and noticing that $(m+k)T \leq t < (m+k+1)T$, we obtain from (5.26) for $\beta := \max\limits_{t \in [0,T]}\left(\int_0^t p(s)\,ds\right)$:

$$\left|\varphi(t,t_0,x_0) - x^*(t)\right| \leq \left|\varphi\big((m+k)T,t_0,x_0\big) - x^*(0)\right| \exp\left(\int_{(m+k)T}^{t} p(s)\,ds\right)$$
$$\leq \left|\varphi(mT,t_0,x_0) - x^*(0)\right| \exp(-kc) \exp\left(\int_0^{t-(m+k)T} p(s+(m+k)T)\,ds\right)$$
$$= \left|\varphi(mT,t_0,x_0) - x^*(0)\right| \exp(-kc) \exp\left(\int_0^{t-(m+k)T} p(s)\,ds\right) \tag{5.27}$$
$$\leq \left|\varphi(mT,t_0,x_0) - x^*(0)\right| \exp(\beta - kc)$$
$$= \left|\varphi(mT,t_0,x_0) - x^*(0)\right| \exp\left(\beta - \left(\left\lfloor \dfrac{t}{T} \right\rfloor - m\right)c\right)$$
$$\leq \left|\varphi(mT,t_0,x_0) - x^*(0)\right| \exp\left(\beta + (m+1)c - c\dfrac{t}{T}\right)$$

Estimate (5.27) shows that $\lim\limits_{t \to +\infty}\left(\left|\varphi(t,t_0,x_0) - x^*(t)\right|\right) = 0$. The proof is complete. ◁

**Proof of Theorem 3:** Since the solution $(x(t), y(t))$ of the initial value problem (3.1), (3.2) with initial condition $(x(0), y(0)) = (x_0, y_0)$ is bounded, it follows that $\omega(x_0, y_0)$ is compact and invariant.

Let $x^*(t)$ be the $T$-periodic solution of (3.1) with $\gamma = \{x^*(t) : t \in [0,T]\} \subseteq \omega(x_0)$. Clearly, there exists a sequence of times $\{t_k \geq 0 : k = 1,2,...\}$ with $\lim(t_k) = +\infty$ and $\lim(x(t_k)) = x^*(0)$. The sequence $\{(x(t_k), y(t_k)) \in \mathbb{R}^n \times \mathbb{R} : k = 1,2,...\}$ is bounded. Consequently, the Bolzano–Weierstrass theorem implies the existence of a subsequence $\{\tau_m \geq 0 : m = 1,2,...\} \subseteq \{t_k \geq 0 : k = 1,2,...\}$ for which $\lim((x(\tau_m), x(\tau_m))) = (x^*(0), y^*) \in \omega(x_0, y_0)$.

Since $\omega(x_0, y_0)$ is invariant, it follows that the solution $(\bar{x}(t), \bar{y}(t))$ of (3.1), (3.2) with initial condition $(\bar{x}(0), \bar{y}(0)) = (x^*(0), y^*)$ is in $\omega(x_0, y_0)$. Consequently, $\bar{x}(t) = x^*(t)$ for all $t \geq 0$ while (by virtue of (3.2)) $\bar{y}(t)$ is a solution of

$$\frac{d\bar{y}}{dt} = g(x^*(t), \bar{y}(t)) \tag{5.28}$$

Since $\omega(x_0, y_0)$ is bounded, it follows that the solution of (5.28) with initial condition $\bar{y}(0) = y^*$ is bounded. By virtue of Theorem 1, there exists a $T$-periodic solution $z^*(t)$ of (5.28) such that $\lim_{t \to +\infty} (\|\bar{y}(t) - z^*(t)\|) = 0$. Notice that $(x^*(t), z^*(t))$ is a $T$-periodic solution of (3.1), (3.2) for which $\lim_{t \to +\infty} (\|(\bar{x}(t), \bar{y}(t)) - (x^*(t), z^*(t))\|) = 0$. Therefore, $\omega(x^*(0), y^*) = \bar{\gamma} = \{(x^*(t), z^*(t)) : t \in [0,T]\}$. Since $(x^*(0), y^*) \in \omega(x_0, y_0)$ the fact that $\omega(x_0, y_0)$ is closed and invariant implies that $\omega(x^*(0), y^*) \subseteq \omega(x_0, y_0)$. Consequently, we obtain that $\bar{\gamma} = \{(x^*(t), z^*(t)) : t \in [0,T]\} \subseteq \omega(x_0, y_0)$. The proof is complete. ◁

Next, we provide the proof of Lemma 1.

**Proof of Lemma 1:** Let arbitrary $R > 0$ and $a \in (0,1]$ be given. Define

$$B = 2R + 1 \tag{5.29}$$

Consider the $2\pi$-periodic time-varying transformation

$$x = \Phi(t, w) \tag{5.30}$$

where

$$\Phi(t, w) := w + \varepsilon q(t, w) \tag{5.31}$$

$$q(t, w) := t\bar{h}(w) - \int_0^t \sin(s) h(w + a\sin(s)) ds \tag{5.32}$$

and

$$\bar{h}(w) := \frac{1}{2\pi} \int_0^{2\pi} \sin(t) h(w + a\sin(t)) dt \tag{5.33}$$

The function $[-B, B] \ni w \mapsto \Phi(t, w)$ is increasing for each $t \geq 0$ since the following inequalities hold for all $t \in [0, 2\pi]$, $a \in (0,1]$ and $w \in [-B, B]$:

$$\frac{\partial \Phi}{\partial w}(t, w) = 1 + \frac{\varepsilon t}{2\pi} \int_0^{2\pi} \sin(s) h'(w + a\sin(s)) ds - \varepsilon \int_0^t \sin(s) h'(w + a\sin(s)) ds$$

$$\geq 1 - 2\varepsilon \int_0^{2\pi} |\sin(s)| |h'(w + a\sin(s))| ds \geq 1 - 2\varepsilon \max_{|z| \leq B+a} (|h'(z)|) \int_0^{2\pi} |\sin(s)| ds \quad (5.34)$$

$$\geq 1 - 8\varepsilon \max_{|z| \leq B+1} (|h'(z)|) > 0$$

The last inequality above is a consequence of (4.7) and (5.29).

Moreover, for all $t \in [0, 2\pi]$ it holds that

$$[-R, R] \subseteq [\Phi(t, -2R), \Phi(t, 2R)] \quad (5.35)$$

Indeed, due to (5.31), (5.32), (5.33) we obtain for all $t \in [0, 2\pi]$ and $a \in (0,1]$:

$$\Phi(t, 2R) = 2R + \frac{\varepsilon t}{2\pi} \int_0^{2\pi} \sin(s) h(2R + a\sin(s)) ds - \varepsilon \int_0^t \sin(s) h(2R + a\sin(s)) ds$$

$$\geq 2R - 8\varepsilon \max_{|z| \leq 2R+1} (|h(z)|) \geq R \quad (5.36)$$

The last inequality above is a consequence of (4.8). Similarly, we obtain $-R \geq \Phi(t, -2R)$ for all $t \in [0, 2\pi]$.

Definition (5.33) implies the following equations for all $w \in \mathbb{R}$:

$$\bar{h}(w) = \frac{1}{2\pi} \int_0^{2\pi} \sin(t) h(w + a\sin(t)) dt = \frac{1}{2\pi} \int_0^{2\pi} \sin(t) (h(w + a\sin(t)) - h(w)) dt$$

$$= \frac{a}{2\pi} \int_0^{2\pi} \int_0^1 \sin^2(t) h'(w + a\lambda \sin(t)) d\lambda dt$$

$$= \frac{ah'(w)}{2} + \frac{a}{2\pi} \int_0^{2\pi} \int_0^1 \sin^2(t) (h'(w + a\lambda \sin(t)) - h'(w)) d\lambda dt$$

$$= \frac{ah'(w)}{2} + \frac{a^2}{2\pi} \int_0^{2\pi} \int_0^1 \int_0^1 \lambda \sin^3(t) h''(w + a\lambda\mu \sin(t)) d\mu d\lambda dt$$

$$= \frac{ah'(w)}{2} + \frac{a^2}{2\pi} \int_0^{2\pi} \int_0^1 \int_0^1 \lambda \sin^3(t) (h''(w + a\lambda\mu \sin(t)) - h''(w)) d\mu d\lambda dt$$

$$= \frac{ah'(w)}{2} + \frac{a^3}{2\pi} \int_0^{2\pi} \int_0^1 \int_0^1 \int_0^1 \lambda^2 \mu \sin^4(t) h'''(w + as\lambda\mu \sin(t)) ds d\mu d\lambda dt$$

Consequently, we obtain for all $w \in [-B, B]$ and $a \in (0,1]$:

$$w\bar{h}(w) \geq \frac{a}{2} wh'(w) - \frac{a^3}{12}|w| \max_{|s| \leq B+1}(|h'''(s)|) \qquad (5.37)$$

Definition (4.10) guarantees that the following implication holds:

$$|w| \leq B \Rightarrow wh'(w) \geq \sigma w^2 \qquad (5.38)$$

Combining (5.37) and (5.38) we get for all $w \in [-B, B]$, $\varepsilon > 0$ and $a \in (0,1]$:

$$-\varepsilon w \bar{h}(w) \leq -\frac{\varepsilon a}{2} \sigma w^2 + \frac{\varepsilon a^3}{12} |w| \max_{|s| \leq B+1}(|h'''(s)|) \qquad (5.39)$$

Let arbitrary $t_0 \geq 0$, $x_0 \in \mathbb{R}$ with $|x_0| \leq R$ be given. Clearly, there exists $t_{\max} \in (t_0, +\infty]$ for which the unique solution $x(t)$ of (4.1) with initial condition $x(t_0) = x_0$ is defined for $t \in [t_0, t_{\max})$. Moreover, if $t_{\max} < +\infty$ then $\limsup_{t \to t_{\max}^-}(|x(t)|) = +\infty$.

Define $w(t)$ so that $x(t) = \Phi(t, w(t))$ for all $t \in [t_0, t_{\max})$ for which $x(t) \in [\Phi(t, -B), \Phi(t, B)]$. Notice that (5.35) and the fact that $|x_0| \leq R$ implies that $|w(t_0)| \leq 2R < B$. By continuity of the solution there exists $\bar{T} > t_0$ such that $x(t) \in [\Phi(t, -B), \Phi(t, B)]$ (or equivalently $|w(t)| \leq B$) for all $t \in [t_0, \bar{T}]$. Using the equation $x(t) = \Phi(t, w(t))$ and equations (4.1), (5.31), (5.32), (5.33) we obtain for all $t \in [t_0, \bar{T}]$:

$$\begin{aligned}
&\frac{\partial \Phi}{\partial w}(t, w(t)) \dot{w}(t) + \varepsilon \bar{h}(w(t)) \\
&= -\varepsilon \sin(t) \Big( h\big(w(t) + \varepsilon q(t, w(t)) + a \sin(t)\big) - h\big(w(t) + a \sin(t)\big) \Big) \\
&= -\varepsilon^2 q(t, w(t)) \sin(t) \int_0^1 h'\big(w(t) + s\varepsilon q(t, w(t)) + a \sin(t)\big) ds
\end{aligned} \qquad (5.40)$$

Since $|w(t)| \leq B$ for all $t \in [t_0, \bar{T}]$, we obtain from (5.32), (5.33) and the fact $a \in (0,1]$ for all $t \in [t_0, \bar{T}]$:

$$\begin{aligned}
|q(t, w(t))| &= \left| \frac{t}{2\pi} \int_0^{2\pi} \sin(s) h(w(t) + a \sin(s)) ds - \int_0^t \sin(s) h(w(t) + a \sin(s)) ds \right| \\
&\leq 2 \int_0^{2\pi} |\sin(s)||h(w(t) + a \sin(s))| ds \leq 2 \max_{|z| \leq B+1}(|h(z)|) \int_0^{2\pi} |\sin(s)| ds \leq 8 \max_{|z| \leq B+1}(|h(z)|)
\end{aligned} \qquad (5.41)$$

Since $|w(t)| \leq B$ for all $t \in [t_0, \bar{T}]$, we obtain from (5.39), (5.40), (5.31) and the fact $a \in (0,1]$ for all $t \in [t_0, \bar{T}]$:

$$\frac{\partial \Phi}{\partial w}(t, w(t)) w(t) \dot{w}(t) \leq -\frac{\varepsilon a}{2} \sigma w^2(t)$$
$$+ \frac{\varepsilon a^3}{12} |w(t)| \max_{|s| \leq B+1} (|h'''(s)|) + 8\varepsilon^2 |w(t)| \max_{|z| \leq B+1} (|h(z)|) \max_{|z| \leq B+R+1} (|h'(z)|) \quad (5.42)$$

Inequalities (5.34), (5.42) guarantee that the following implication is valid for all $t \in [t_0, \bar{T}]$:

$$|w(t)| > \frac{a^2}{6\sigma} \max_{|s| \leq B+1} (|h'''(s)|) + \frac{16\varepsilon}{a\sigma} \max_{|z| \leq B+1} (|h(z)|) \max_{|z| \leq B+R+1} (|h'(z)|) \Rightarrow w(t)\dot{w}(t) < 0 \quad (5.43)$$

The fact that $\frac{a^2}{6\sigma} \max_{|s| \leq B+1} (|h'''(s)|) + \frac{16\varepsilon}{a\sigma} \max_{|z| \leq B+1} (|h(z)|) \max_{|z| \leq B+R+1} (|h'(z)|) \leq R$ (a consequence of (4.9) and (5.29)) and implication (5.43) guarantees that $|w(t)| \leq \max(R, |w(t_0)|) \leq 2R < B$ for all $t \in [t_0, \bar{T}]$. A standard contradiction argument allows us to conclude that $|w(t)| \leq \max(R, |w(t_0)|) \leq 2R < B$ for all $t \geq t_0$. This implies that $x(t) \in [\Phi(t, -2R), \Phi(t, 2R)]$ for all $t \geq t_0$ and using the fact that $x(t) = \Phi(t, w(t))$, (5.31) and (5.41) we get for all $t \geq t_0$:

$$|x(t)| = |w(t) + \varepsilon q(t, w(t))|$$
$$\leq |w(t)| + \varepsilon |q(t, w(t))| \leq |w(t)| + 8\varepsilon \max_{|z| \leq B+1} (|h(z)|) \quad (5.44)$$

Implication (5.43) guarantees that

$$\limsup_{t \to +\infty} (|w(t)|) \leq \frac{a^2}{6\sigma} \max_{|s| \leq B+1} (|h'''(s)|) + \frac{16\varepsilon}{a\sigma} \max_{|z| \leq B+1} (|h(z)|) \max_{|z| \leq B+R+1} (|h'(z)|) \quad (5.45)$$

Estimates (5.44), (5.45) and definition (5.29) guarantee estimate (4.11).

The proof is complete. ◁

We finish this section with the proof of Theorem 3 and the proof of Proposition 1.

**Proof of Theorem 4:** Inequalities (4.12), (4.13) guarantee the existence of $\varepsilon^* > 0$ and $b > 0$ such that (4.5), (4.7), (4.8), (4.9) hold for all $\varepsilon \in (0, \varepsilon^*)$ and

$$\frac{a^2 \gamma}{6\sigma} + \varepsilon K < b \quad (5.46)$$

The rest of the proof is a consequence of Theorem 1, Theorem 2 (with $S = [-R, R] \subseteq \mathbb{R}$ and $C = [-b, b]$), Lemma 1 (which guarantees that $\limsup_{t \to +\infty} (|x(t)|) < b$) and definitions $\gamma := \max_{|z| \leq 2R+2} (|h'''(z)|)$, $K := 8 \max_{|z| \leq 2R+2} (|h(z)|) \left( \frac{2}{a\sigma} \max_{|z| \leq 3R+2} (|h'(z)|) + 1 \right)$. The proof is complete. ◁

**Proof of Proposition 1:** Define the closed subset of the Banach space $C^0([0, \pi])$ with norm $\|x\|_\infty = \max_{t \in [0, \pi]} (|x(t)|)$ for $x \in C^0([0, \pi])$:

$$B_R = \left\{ x \in C^0([0, \pi]) : \|x\|_\infty = \max_{t \in [0, \pi]} (|x(t)|) \leq R \right\} \tag{5.47}$$

Define the map $P : C^0([0, \pi]) \to C^0([0, \pi])$ by means of the formula for $x \in C^0([0, \pi])$ and $t \in [0, \pi]$:

$$(P(x))(t) := \frac{\varepsilon}{2} \int_t^\pi \sin(s) h(x(s) + a \sin(s)) ds - \frac{\varepsilon}{2} \int_0^t \sin(s) h(x(s) + a \sin(s)) ds \tag{5.48}$$

Definitions (5.47), (5.48) show that the following estimate holds for all $x \in B_R$:

$$\|P(x)\|_\infty \leq \frac{\varepsilon \pi}{2} \max_{|z| \leq R+a} (|h(z)|) \tag{5.49}$$

It follows from (5.49) and (4.14) that $P(B_R) \subseteq B_R$.

We also get from (5.48) for all $x, y \in C^0([0, \pi])$ and $t \in [0, \pi]$:

$$|(P(x))(t) - (P(y))(t)| \leq \frac{\varepsilon}{2} \int_t^\pi \sin(s) |h(x(s) + a \sin(s)) - h(y(s) + a \sin(s))| ds$$
$$+ \frac{\varepsilon}{2} \int_0^t \sin(s) |h(x(s) + a \sin(s)) - h(y(s) + a \sin(s))| ds \tag{5.50}$$

Exploiting (5.47), (5.50) we obtain for all $x, y \in B_R$:

$$\|P(x) - P(y)\|_\infty \leq \frac{\varepsilon \pi}{2} \max_{|z| \leq R+a} (|h'(z)|) \|x - y\|_\infty \tag{5.51}$$

It follows from (5.51) and (4.14) that the map $P : B_R \to B_R$ is a contraction. Banach's fixed point theorem implies the existence of $x^{**} \in B_R$ such that $P(x^{**}) = x^{**}$. Equivalently, this means that the following equation holds for all $t \in [0, \pi]$:

$$x^{**}(t) = \frac{\varepsilon}{2}\int_t^\pi \sin(s)h\big(x^{**}(s)+a\sin(s)\big)ds - \frac{\varepsilon}{2}\int_0^t \sin(s)h\big(x^{**}(s)+a\sin(s)\big)ds \qquad (5.52)$$

Equation (5.52) shows that $x^{**}(0) = -x^{**}(\pi) > 0$. Exploiting (5.52), we get for all $t \in [0,\pi]$:

$$\dot{x}^{**}(t) = -\varepsilon \sin(t) h\big(x^{**}(t)+a\sin(t)\big) \qquad (5.53)$$

Thus, $x^{**}$ is a solution of (4.1) on $[0,\pi]$. We next make the following extension of $x^{**}$ on $[\pi, 2\pi]$ by defining for all $t \in (0,\pi]$:

$$x^{**}(t+\pi) = -x^{**}(t) \qquad (5.54)$$

Using the fact that $h$ is an even function, we obtain from (5.53), (5.54) for all $t \in [\pi, 2\pi]$:

$$\begin{aligned}\dot{x}^{**}(t+\pi) &= -\dot{x}^{**}(t) = \varepsilon \sin(t)h\big(x^{**}(t)+a\sin(t)\big) \\ &= -\varepsilon \sin(t+\pi)h\big(-x^{**}(t+\pi)-a\sin(t+\pi)\big) \\ &= -\varepsilon \sin(t+\pi)h\big(x^{**}(t+\pi)+a\sin(t+\pi)\big)\end{aligned}$$

Thus the $2\pi$ – periodic extension of $x^{**}$ is a $2\pi$ – periodic solution of (4.1). The fact that $\max_{t\in[0,2\pi]}\big(|x^{**}(t)|\big) \leq \frac{\varepsilon\pi}{2}\max_{|z|\leq R+a}\big(|h(z)|\big)$ is a direct consequence of (5.49) and (5.54). The proof is complete. ◁

## 6. Conclusions and Future Work

The paper provided two results that are useful in the study of the existence and the stability properties of a periodic solutions. Both results were applied to extremum seeking problems for a static output map with no plant dynamics and novel non-local results were provided without the use of averaging theorems and singular perturbation arguments.

It is clear that additional results are needed for the study of more general extremum seeking problems: here we focused only on the case of a static map of one variable with no plant dynamics. Moreover, different architectures for the extremum seeking control (e.g., higher order extremum seeking schemes) cannot be covered by the results of the present work. This will be the topic of future research.